\documentclass[11pt]{amsart}
\usepackage{amssymb}

\newcommand{\Sym}{{\mathfrak{S}}}
\newcommand{\Z}{{\mathbb{Z}}}
\newcommand{\Q}{{\mathbb{Q}}}
\newcommand{\F}{{\mathbb{F}}}
\newcommand{\N}{{\mathbb{N}}}

\newcommand{\C}{{\mathbb{C}}}

\newcommand{\ba}{{\mathbf{a}}}
\newcommand{\bc}{{\mathbf{c}}}
\newcommand{\bu}{{\mathbf{u}}}

\newcommand{\cB}{{\mathcal{B}}}

\newcommand{\cM}{{\mathcal{M}}}
\newcommand{\cK}{{\mathcal{K}}}
\newcommand{\bH}{{\mathbf{H}}}
\newcommand{\bT}{{\mathbf{T}}}
\newcommand{\cF}{{\mathfrak{F}}}
\newcommand{\cU}{{\mathcal{U}}}

\newcommand{\fsl}{{\mathfrak{sl}}}

\newcommand{\domi}{\trianglelefteq}

\newcommand\Ind{\operatorname{Ind}}
\newcommand{\Irr}{{\operatorname{Irr}}}

\newcommand{\ulambda}{{\boldsymbol{\lambda}}}
\newcommand{\umu}{{\boldsymbol{\mu}}}
\newcommand{\unu}{{\boldsymbol{\nu}}}
\newcommand{\uvar}{{\underline{\varnothing}}}

\renewcommand{\leq}{\leqslant}
\renewcommand{\geq}{\geqslant}
\renewcommand{\atop}[2]{\genfrac{}{}{0pt}{}{#1}{#2}}

\newtheorem{thm}{Theorem}[section]
\newtheorem{lem}[thm]{Lemma}
\newtheorem{cor}[thm]{Corollary}
\newtheorem{prop}[thm]{Proposition}

\theoremstyle{definition}
\newtheorem{exmp}[thm]{Example}
\newtheorem{defn}[thm]{Definition}

\theoremstyle{remark}
\newtheorem{rem}[thm]{Remark}

\begin{document}

\title{Canonical basic sets in type $B_n$}
\author{Meinolf Geck and Nicolas Jacon}
\address{M.G.: Department of Mathematical Sciences, King's College,
Aberdeen University, Aberdeen AB24 3UE, Scotland, U.K.}

\email{m.geck@maths.abdn.ac.uk}

\address{N.J.: Universit\'e de Franche-Comt\'e,  UFR Sciences et
Techniques, 16 route de Gray, 25 030 Besan\c{c}on, France.}

\email{jacon@math.univ-fcomte.fr}

\date{January, 2006}

\dedicatory{To Gordon James on his $60$th birthday}

\begin{abstract}
More than $10$ years ago, Dipper, James and Murphy developped the theory
of Specht modules for Hecke algebras of type $B_n$. More recently, using
Lusztig's $\ba$-function, Geck and Rouquier showed how to obtain
parametrisations of the irreducible representations of Hecke algebras (of
any finite type) in terms of so-called {\em canonical basic sets}. For
certain values of the parameters in type $B_n$, combinatorial descriptions
of these basic sets were found by Jacon, based on work of Ariki and
Foda--Leclerc--Okado--Thibon--Welsh. Here, we consider the canonical basic
sets for all the remaining choices of the parameters.
\end{abstract}

\maketitle

\pagestyle{myheadings}

\markboth{Geck and Jacon}{Canonical basic sets in type $B_n$}

\maketitle

\section{Introduction: Specht modules in type $B_n$} \label{sec:intro}
Let $W_n$ be the finite Weyl group of type $B_n$, with generating set
$S_n=\{t, s_1,\ldots,s_{n-1}\}$ and relations given by the following diagram:
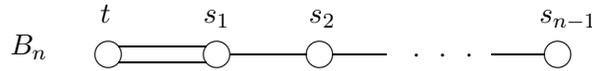
\begin{center}
\begin{picture}(250,30)
\put(  3, 08){$B_n$}
\put( 40, 08){\circle{10}}
\put( 44, 05){\line(1,0){33}}
\put( 44, 11){\line(1,0){33}}
\put( 81, 08){\circle{10}}
\put( 86, 08){\line(1,0){29}}
\put(120, 08){\circle{10}}
\put(125, 08){\line(1,0){20}}
\put(155, 05){$\cdot$}
\put(165, 05){$\cdot$}
\put(175, 05){$\cdot$}
\put(185, 08){\line(1,0){20}}
\put(210, 08){\circle{10}}
\put( 37, 20){$t$}
\put( 76, 20){$s_1$}
\put(116, 20){$s_2$}
\put(203, 20){$s_{n-1}$}
\end{picture}
\end{center}
We have $W_n \cong ({\Z}/2\Z)^n \rtimes \Sym_n$, where $\Sym_n$ is the
symmetric group on $n$ letters. Let $k$ be a field and $Q,q\in k^\times$.
We denote by $H_n:=H_k(W_n,Q,q)$ the corresponding Iwahori--Hecke algebra.
This is an associative algebra over $k$, with a basis $\{T_w\mid w \in W_n\}$
such that the following relations hold for the multiplication:
\begin{align*}
T_t^2 & = Q\, T_1+(Q-1)T_t,\\
T_{s_i}^2 &= q\, T_1+ (q-1)T_{s_i} \qquad \mbox{for $1 \leq i \leq n-1$},\\
T_w & = T_{r_1} \cdots T_{r_m}, \qquad \mbox{for $w\in W_n$, $r_i \in S_n$,
$m=l(w)$},
\end{align*}
where $\ell\colon W_n \rightarrow \N$ is the usual length function.
(We write $\N=\{0,1,2,\ldots\}$.) In this paper, we are concerned with
the problem of classifying the irreducible representations of $H_n$.
For applications to the representation theory of finite groups of Lie
type, see Lusztig's book \cite{Lu4} and the surveys \cite{DGHM},
\cite{mylaus}.

Hoefsmit \cite{Hoefs} explicitly constructed the irreducible
representations of $H_n$ in the case where $H_n$ is semisimple. In
this case, we have a parametrization
\[ \Irr(H_n)=\{ E^\ulambda \mid \ulambda \in \Pi_n^2\}.\]
Here, $\Pi_n^2$ denotes the set of all bipartitions of $n$, that is,
pairs of partitions $\ulambda=(\lambda^{(1)}\,|\,\lambda^{(2)})$ such that
$|\lambda^{(1)}|+|\lambda^{(2)}|=n$. In the general case where $H_n$ is
not necessarily semisimple, Dipper--James--Murphy \cite{DJM3} constructed
a {\em Specht module}\footnote{Actually, throughout this paper, we
shall denote by $S^{\ulambda}$ the module that is labelled by
$({\lambda^{(2)}}^*, {\lambda^{(1)}}^*)$ in \cite{DJM3}, where the star
denotes the conjugate partition. For example, the index representation
(given by $T_t\mapsto Q$, $T_{s_i}\mapsto q$) is labelled by the pair
$((n), \varnothing)$ and the sign representation (given by $T_t
\mapsto -1$, $T_{s_i} \mapsto -1$) is labelled by $(\varnothing,(1^n))$.
Thus, our labelling coincides with that of Hoefsmit \cite{Hoefs}.}
$S^\ulambda$ for any bipartition $\ulambda\in \Pi_n^2$. This is an
$H_n$-module, reducible in general, such that $\dim S^\ulambda=
\dim E^\ulambda$.
Every Specht module comes equipped with an $H_n$-equivariant symmetric
bilinear form. Taking the quotient by the radical, we obtain an $H_n$-module
\[ D^{\ulambda}:=S^{\ulambda}/ \mbox{rad}(S^{\ulambda}) \qquad \mbox{for
any $\ulambda \in \Pi_n^2$}.\]
Let $\Lambda_n^2:=\{\ulambda\in \Pi_n^2\mid D^\ulambda\neq
\{0\}\}$.  Then, by \cite[\S 6]{DJM3}, we have
\[\Irr(H_n)=\{D^{\ulambda}\mid \ulambda \in \Lambda_n^2\}.\]
Through the work of Dipper--James--Murphy \cite{DJM3}, Ariki--Mathas
\cite{ArMa} and Ariki \cite{Ar2}, \cite{Ar3}, we have an explicit
description of the set $\Lambda_n^2$. This depends on the
following two parameters:
\begin{align*}
e&:=\min\{i\geq 2\mid 1+q+q^2+\cdots +q^{i-1}=0\}\\
f_n(Q,q)&:= \prod_{i=-(n-1)}^{n-1} (Q+q^i).
\end{align*}
(We set $e=\infty$ if $1+q+\cdots +q^{i-1} \neq 0$ for all $i\geq 2$.)

\begin{thm} \label{mainspecht} The set $\Lambda_n^2$ is given as
follows.
\begin{itemize}
\item[(A)] {\rm (Dipper--James \cite[4.17 and 5.3]{DJ2} and
Dipper--James--Murphy \cite[6.9]{DJM3}).} Assume that $f_n(Q,q) \neq 0$. Then
\[ \Lambda_n^2=\{\ulambda\in \Pi_n^2\mid \mbox{$\lambda^{(1)}$
and $\lambda^{(2)}$ are $e$-regular}\}.\]
(A partition is called $e$-regular if no part is repeated $e$ times or
more.)
\item[(B)] {\rm (Dipper--James--Murphy \cite[7.3]{DJM3}).}
Assume that $f_n(Q,q)=0$ and $q=1$ (and, hence, $Q=-1$). Then
\[ \Lambda_n^2=\{\ulambda \in \Pi_n^2\mid \mbox{$\lambda^{(1)}$
is $e$-regular and $\lambda^{(2)}=\varnothing$}\}.\]
\item[(C)] {\rm (Ariki--Mathas \cite{ArMa}, Ariki \cite{Ar2}).} Assume
that $f_n(Q,q)=0$ and $q\neq 1$. Thus, $Q=-q^d$ where $-(n-1)\leq d
\leq n-1$. Then
\[ \Lambda_n^2=\{ \ulambda \in \Pi_n^2\mid \ulambda \mbox{ is a
{\em Kleshchev} $e$-bipartition}\};\]
this set only depends on $e$ and $d$. (See Remark~\ref{kleshdef} where we
recall the exact definition of Kleshchev bipartitions.)
\end{itemize}
\end{thm}

Now, a fundamental feature of an Iwahori--Hecke algebra as above is that
it can be derived from a ``generic'' algebra by a process of
``specialisation''. For this purpose, let us assume that we can write
\begin{equation*}
\mbox{$Q=\xi^b$ and $q=\xi^a$, where $\xi \in k^\times$ and $a,b \in \N$}.
\tag{$\clubsuit$}
\end{equation*}
Thus, the parameters $Q,q$ of $H_n$ are assumed to be integral powers of
one fixed non-zero element $\xi$ of~$k$. (Note that $\xi,a,b$ are not
uniquely determined by $Q,q$.) This situation naturally occurs, for example,
in applications to the representation theory of reductive groups over a
finite field $\F_q$, where $\xi=q1_k$. The integers $a,b$ uniquely define
a weight function $L \colon W_n \rightarrow \Z$ in the sense of Lusztig
\cite{Lusztig03}, that is, we have
\begin{gather*}
L(t)=b, \qquad L(s_1)=\cdots =L(s_{n-1})=a,\\
L(ww')=L(w)+L(w') \qquad \mbox{whenever $\ell(ww')= \ell(w)+\ell(w')$}.
\end{gather*}
Let $R={\Z}[u,u^{-1}]$ be the ring of Laurent polynomials in an
indeterminate~$u$. Let $\bH_n=\bH_R(W_n,L)$ be the corresponding
generic Iwahori--Hecke algebra. This is an associative algebra over $R$,
with a free $R$-basis $\{\bT_w\mid w \in W_n\}$ such that the following
relations hold for the multiplication:
\begin{align*}
\bT_t^2 & = u^b\, \bT_1+(u^b-1)\bT_t,\\
\bT_{s_i}^2 &= u^a\, \bT_1+ (u^a-1)\bT_{s_i} \qquad
\mbox{for $1 \leq i \leq n-1$},\\
\bT_w & = \bT_{r_1} \cdots \bT_{r_m}, \qquad
\mbox{for $w\in W_n$, $r_i \in S_n$, $m=l(w)$}.
\end{align*}
Then there is a ring homomorphism $\theta \colon R \rightarrow k$ such that
$\theta(u)=\xi$, and we obtain $H_n$ by extension of scalars from $R$
to $k$ via $\theta$:
\[ H_n=k \otimes_R \bH_R(W_n,L).\]
Now, it would be desirable to obtain a parametrization of the
irreducible representations of $H_n$ which also takes into account the
weight function $L$. (Note that $\Lambda_n^2$ is ``insensitive'' to $L$:
it does not depend on the choice of $a,b,\xi$ such that ($\clubsuit$)
holds.) It would also be desirable to obtain a parametrization which fits
into a general framework valid for Iwahori--Hecke algebras of any finite
type. (Note that it is not clear how to define Specht modules for algebras
of exceptional type, for example.) Such a general framework for obtaining
$L$-adapted parametrizations was developped by Geck \cite{mykl}, \cite{my00},
\cite{mylaus} and Geck--Rouquier \cite{GeRo2}. It relies on deep (and
conjectural for general choices of $L$) properties of the
{\em Kazhdan--Lusztig basis} and Lusztig's $\ba$-{\em function}. In this
framework, the parametrization is in terms of so-called ``canonical basic
sets''. We recall the basic ingredients in Section~\ref{sec:afunc}.

Jacon \cite{Jac0}, \cite{Jac1}, \cite{Jac2} explicitly described these
canonical basic sets in type $B_n$ for certain choices of $a$ and $b$,
most notably the case where $a=b$ (the ``equal parameter case'') and the
case where $b=0$ (which gives a classification of the irreducible
representations of an algebra of type $D_n$). Note that, in these cases,
the canonical basic sets are different from the set $\Lambda_n^2$.
It is known, see Theorem~\ref{asymp}, that $\Lambda_n^2$ can be
interpreted as a canonical basic set with respect to weight functions $L$
such that $b>(n-1)a>0$. The aim of this article to determine the canonical
basic sets for all the remaining choices of $a,b$. This goal
will be achieved in Theorems~\ref{caseA}, \ref{caseB} and \ref{caseC};
note, however, that for ground fields of positive characteristic, our
solution in case~(C) of Theorem~\ref{mainspecht} relies on the validity of
Lusztig's conjectures \cite{Lusztig03} on Hecke algebras with unequal
parameters.

An application of the results in this paper to the modular representation
theory of finite groups of Lie type can be found in \cite{myprinc}: The 
explicit description of canonical basic sets for $H_n$ yields a natural
parametrization of the modular principal series representations for 
finite classical groups.

\section{Lusztig's $\ba$-function and the decomposition matrix}
\label{sec:afunc}

We keep the setting of the previous section, where $\bH_n=\bH_R(W_n,L)$
is the generic Iwahori--Hecke algebra corresponding to the Weyl group
$W_n$ and the weight function $L$ such that $L(t)=b\geq 0$ and $L(s_i)=a
\geq 0$ for $1\leq i \leq n-1$. The aim of this section is to recall the
basic ingredients in the definition of a {\em canonical basic set} for
the algebra $H_n=k \otimes_R \bH_n$ where $\theta \colon R \rightarrow k$
is a ring homomorphism into a field $k$ and $\xi=\theta(u)$.

At the end of this section, in Theorem~\ref{asymp}, we recognise the set
$\Lambda_n^2$ (arising from the theory of Specht modules) as a canonical
basic set for a certain class of weight functions.

Let $K={\Q}(u)$ be the field of fractions of $R$. By extension of
scalars, we obtain a $K$-algebra $\bH_{K,n}=K \otimes_R \bH_n$, which is
known to be split semisimple; see Dipper--James \cite{DJ2}. Furthermore,
in this case, we have
\[ \Irr(\bH_{K,n})=\{S^\ulambda \mid \ulambda \in \Pi_n^2\}\]
where $S^\ulambda$ are the Specht modules defined by Dipper--James--Murphy
\cite{DJM3}. (Recall our convention about the labelling of these modules.)

The definition of Lusztig's $\ba$-function relies on the fact that
$\bH_n$ is a symmetric algebra. Indeed, we have a trace form $\tau
\colon \bH_n \rightarrow R$ defined by
\[ \tau(\bT_1)=1 \qquad \mbox{and}\qquad \tau(\bT_w)=0 \quad \mbox{for
$w\neq 1$}.\]
The associated bilinear form $\bH_n \times \bH_n \rightarrow R$,
$(h,h')\mapsto \tau(hh')$, is symmetric and non-degenerate.
Thus, $\bH_n$ has the structure of a {\em symmetric algebra}.
Now extend $\tau$ to trace form $\tau_K \colon \bH_{K,n} \rightarrow K$.
Since $\bH_{K,n}$ is split semisimple, we have
\[ \tau_K(T_w)=\sum_{\ulambda \in \Pi_n^2} \frac{1}{\bc_\ulambda} \,
\mbox{trace}(\bT_w,S^\ulambda) \quad \mbox{for all $w\in W_n$},\]
where $0 \neq \bc_\ulambda\in R$; see \cite[Chapter~7]{ourbuch}.
We have
\[ \bc_\ulambda=f_\ulambda u^{-\ba_\ulambda} + \mbox{combination of higher
powers of $u$},\]
where both $f_\ulambda$ and $\ba_\ulambda$ are integers, $f_\ulambda>0$,
$\ba_\ulambda\geq 0$.

Hoefsmit \cite{Hoefs} obtained explicit combinatorial formulas for
$\bc_\ulambda$. Then Lusztig deduced purely combinatorial expressions
for $f_\ulambda$ and $\ba_\ulambda$; see \cite[Chap.~22]{Lusztig03}. Let
us first assume that $a>0$.  Then we have
\begin{alignat*}{2}
f_\ulambda  &=1 \quad &&\mbox{if $a>0$ and $b/a\not\in\{0,1,\ldots,n-1\}$},\\
f_\ulambda & \in \{1,2,4,8,16,\ldots\} \quad &&\mbox{if $a>0$ and $b/a
\in \{0,1,\ldots,n-1\}$};
\end{alignat*}
see \cite[22.14]{Lusztig03}. To describe $\ba_\ulambda$, we need some
more notation. Let us write
\[ b=ar+b' \qquad \mbox{where $r,b'\in \N$ and $b'<a$}.\]
Let $\ulambda=(\lambda^{(1)}\,|\,\lambda^{(2)})\in \Pi_n^2$ and write
\[ \lambda^{(1)}=(\lambda^{(1)}_{1}\geq \lambda^{(1)}_{2}\geq
\lambda^{(1)}_{3}\geq \cdots), \qquad \lambda^{(2)}=(\lambda^{(2)}_{1}
\geq \lambda^{(2)}_{2}\geq \lambda^{(2)}_{3}\geq \cdots),\]
where $\lambda^{(1)}_{N}=\lambda^{(2)}_{N}=0$ for all large values of $N$.
Now fix a large $N$ such that $\lambda^{(1)}_{N+r+1}=\lambda^{(2)}_{N+1}=0$.
Then we set
\begin{align*}
\alpha_i&=a(\lambda^{(1)}_{N+r-i+1}+i-1)+b'
\qquad\mbox{for $1\leq i\leq N+r$},\\
\beta_j&=a(\lambda^{(2)}_{N-j+1}+j-1)\qquad\mbox{for $1\leq j\leq N$}.
\end{align*}
We have $0 \leq \alpha_1<\alpha_2<\cdots <\alpha_{N+r}$, $0 \leq \beta_1<
\beta_2<\cdots <\beta_N$.  Since $N$ is large, we have
$\alpha_i=a(i-1)+b'$ for $1\leq i \leq r$ and $\beta_j=a(j-1)$ for $1 \leq
j \leq r$. Now we can state:

\begin{prop}[Lusztig \protect{\cite[22.14]{Lusztig03}}] \label{afuncWn}
Recall that $a>0$ and $b=ar+b'$ as above. Then $\ba_\ulambda=A_N-B_N$
where
\begin{align*}
A_N&=\sum_{\atop{1 \leq i \leq N{+}r}{1\leq j\leq N}} \min(\alpha_i,\beta_j)
+\sum_{1 \leq i<j\leq N{+}r} \min(\alpha_i,\alpha_j)
+\sum_{1 \leq i<j\leq N} \min(\beta_i,\beta_j),\\
B_N&=\sum_{\atop{1 \leq i \leq N{+}r}{1 \leq j \leq N}}
\min(a(i-1)+b',a(j-1))\\
& \qquad +\sum_{1 \leq i<j\leq N{+}r} \min(a(i-1)+b',a(j-1)+b')\\
&\qquad\qquad +\sum_{1 \leq i<j\leq N} \min(a(i-1),a(j-1)).
\end{align*}
(Note that $B_N$ only depends on $a,b,n,N$ but not on $\ulambda$.)
\end{prop}

\begin{rem} \label{expa} (a) Assume that $a=0$ and $b>0$. Then one easily
checks, directly using Hoesfmit's formulas for $\bc_\ulambda$, that
\begin{align*}
f_\ulambda & \in \{\mbox{divisiors of $n!$}\} \qquad \mbox{for all $\ulambda
\in \Pi_n^2$},\\
\ba_\ulambda & =b\,|\lambda^{(2)}|\qquad \mbox{for all $\ulambda
\in \Pi_n^2$}.
\end{align*}
(If $a=b=0$, then $f_\ulambda=|W_n|/\dim E^\ulambda$ and $\ba_\lambda=0$
for all $\ulambda \in \Lambda_n^2$.)

(b) For any $a,b\in \N$, we have
\[ \ba((n), \varnothing)=0 \qquad \mbox{and}\qquad \ba(\varnothing,(1^n))=
L(w_0),\]
where $w_0\in W_n$ is the unique element of maximal length.
\end{rem}

\begin{defn} \label{goodL} We say that $k$ is {\em $L$-good} if
$f_\ulambda 1_k\neq 0$ for all $\ulambda \in \Pi_n^2$. In
particular, any field of characteristic $0$ is $L$-good. For
fields of characteristic $p>0$, the above formulas for $f_\ulambda$
show that the conditions are as follows.
\begin{itemize}
\item[(i)] If $a>0$ and $b/a\not\in\{0,1,\ldots,n-1\}$ then any field is
$L$-good.
\item[(ii)] If $a>0$ and $b/a\in\{0,1,\ldots,n-1\}$, then fields of
characteristic $p\neq 2$ are $L$-good.
\item[(iii)] If $a=0$, then fields of characteristic $p>n$ are $L$-good.
\end{itemize}
(Note that the case $a=0$ should be merely considered as a curiosity,
which may only show up in extremal situations as far as applications
are concerned.)
\end{defn}

Now let us consider the {\em decomposition matrix} of $H_n$,
\[ D=\bigl([S^\ulambda:D^\umu]\bigr)_{\ulambda \in \Pi_n^2,\umu \in
\Lambda_n^2},\]
where $[S^\ulambda:D^\umu]$ denotes the multiplicity of $D^\umu$ as a
composition factor of $S^\ulambda$. By Dipper--James--Murphy
\cite[\S 6]{DJM3}, we have
\begin{equation*}
\left\{\begin{array}{l} \quad [S^\umu:D^\umu]=1 \quad \mbox{for any
$\umu\in \Lambda_n^2$},\\ \quad [S^\ulambda:D^\umu] \neq 0 \quad
\Rightarrow \quad \ulambda \trianglelefteq \umu,\end{array}\right.
\tag{$\Delta$}
\end{equation*}
where $\trianglelefteq$ denotes the dominance order on bipartitions. Note
that these conditions uniquely determine the set $\Lambda_n^2$
once the matrix $D$ is known.

\begin{defn} \label{defbs} Let $\beta \colon \Lambda_n^2
\rightarrow \Pi_n^2$ be an injective map and set $\cB:=
\beta(\Lambda_n^2) \subseteq \Pi_n^2$. Let us denote
\[ M^{\unu}:=D^{\beta^{-1}(\unu)} \qquad \mbox{for any $\unu \in \cB$}.\]
We say that $\cB$ is ``canonical basic set'' for $H_n$ if the following
conditions are satisfied:
\begin{equation*}
\left\{\begin{array}{l} \quad [S^{\unu}:M^{\unu}]=1 \quad \mbox{for any
$\unu\in \cB$},\\ \quad [S^{\ulambda}:M^\unu] \neq 0 \quad
\Rightarrow \quad \ulambda=\unu \mbox{ or } \ba_\unu <\ba_\ulambda,
\end{array}\right.  \tag{$\Delta_\ba$}
\end{equation*}
Note that the conditions ($\Delta_\ba$) uniquely determine the set $\cB$
and the bijection $\beta \colon \Lambda_n^2 \stackrel{\sim}{\rightarrow}
\cB$. Thus, we obtain a new ``canonical'' labelling
\[ \Irr(H_n)=\{M^\unu \mid \nu \in \cB\}.\]
Furthermore, the submatrix
\[ D^\circ=\bigl([S^\unu:M^{\unu'}]\bigr)_{\unu,\unu' \in \cB}\]
is square and lower triangular with $1$ on the diagonal, when we order
the modules according to increasing values of $\ba_\unu$. More precisely,
we have a block lower triangular shape
\[ D^\circ=\left(\begin{array}{cccc} D_0^\circ &&&0 \\
& D_1^{\circ} &&\\ && \ddots & \\ *&&& D_N^\circ \end{array}\right),\]
where the block $D_i^\circ$ has rows and columns labelled by those
$S^\unu$ and $M^{\unu'}$, respectively, where $\ba_\unu=\ba_{\unu'}=i$,
and each $D_i^{\circ}$ is the identity matrix.
\end{defn}

\begin{thm}[Geck \cite{mykl}, \cite{my00}, \protect{\cite[\S 6]{mylaus}} 
and Geck--Rouquier \cite{GeRo2}] \label{gero} Assume that Lusztig's
conjectures {\bf (P1)--(P14)} in \cite[14.2]{Lusztig03} and a certain
weak version of {\bf (P15)} (as specified in \cite[5.2]{mylaus}) hold
for $\bH_n=\bH_A(W_n,L)$. Assume further that $k$ is $L$-good; see
Definition~\ref{goodL}. Then $H_n$ admits a canonical basic set.
\end{thm}

\begin{rem} (a) The above result is proved by a general argument which
works for Iwahori--Hecke algebras of any finite type, once the properties
{\bf (P1)--(P14)} and the weak version of {\bf (P15)} are known to hold.
This is the case, for example, when the weight function $L$ is a multiple
of the length function (the ``equal parameter case''); see Lusztig
\cite[Chap.~15]{Lusztig03}. However, it seems to be very hard to
obtain an explicit description of $\cB$ from the construction in the proof.

(b) If $k$ is not $L$-good, it is easy to produce examples in which
a canonical basic set does not exist; see \cite[4.15]{mylaus}.
\end{rem}

\begin{table}[htbp]
\caption{Decomposition numbers for $B_3$ with $Q=1$, $q=-1$.}
\label{tabb3}
\begin{center}
$\renewcommand{\arraystretch}{1.1}
\begin{array}{|c|cccc|} \hline  &
\multicolumn{4}{c|}{[S^\ulambda:D^\umu]} \\ \hline
(3\,|\,\varnothing)     & 1 & . & . & . \\
(21\,|\,\varnothing)    & . & 1 & . & . \\
(111\,|\,\varnothing)   & 1 & . & . & . \\
(2\,|\,1)               & . & 1 & 1 & . \\
(11\,|\,1)              & . & 1 & 1 & . \\
(1\,|\,2)               & 1 & . & . & 1 \\
(\varnothing\,|\,3 )    & . & . & 1 & . \\
(1\,|\,11)             & 1 & . & . & 1 \\
(\varnothing\,|\,21)   & . & . & . & 1 \\
(\varnothing\,|\,111)  & . & . & 1 & . \\
\hline\end{array} \qquad
\begin{array}{|c|cc|}
\hline \ba_\ulambda& b=0 & b=4 \\ \hline
(3\,|\,\varnothing)    & 0 & 0   \\
(21\,|\,\varnothing)   & 2 & 1   \\
(111\,|\,\varnothing)  & 6 & 3   \\
(2\,|\,1)              & 1 & 4   \\
(11\,|\,1)             & 3 & 5   \\
(1\,|\,2)              & 1 & 7   \\
(\varnothing\,|\,3 )   & 0 & 9   \\
(1\,|\,11)             & 3 & 10  \\
(\varnothing\,|\,21)   & 2 & 13  \\
(\varnothing\,|\,111)  & 6 & 18  \\ \hline \end{array}$
\end{center}
\end{table}

\begin{exmp} \label{exb3} Assume that $n=3$, $Q=1$ and $q=-1$. Thus,
in $H_3$, we have the quadratic relations
\[ T_t^2 = T_1 \quad \mbox{and}\quad T_{s_i}^2 = -T_1-2T_{s_i}
\quad \mbox{for $i=1,2$}.\]
In this case, we have
\[ \Lambda_3^2=\{(3\,|\,\varnothing),\; (21)\,|\,\varnothing),\;
(2\,|\,1),\;(1\,|\,2)\}.\]
The decomposition matrix is printed in Table~\ref{tabb3}.

\begin{itemize}
\item[(a)] Now let us take $\xi=-1$ and write $Q=\xi^b$, $q=\xi^a$ where
$a=1$ and $b=0$. Then the canonical basic set is given by
\[ \cB=\{(3\,|\,\varnothing),\; (\varnothing\,|\,3),\; (1\,|\,2),\;
(2\,|\,1)\}.\]
The bijection $\beta \colon \Lambda_3^2\rightarrow \cB$ is given by
\begin{align*}
(3\,|\,\varnothing) &\mapsto (3\,|\,\varnothing), \\
(2\,|\,\varnothing) & \mapsto (2\,|\,1),\\
(2\,|\,1) & \mapsto (\varnothing\,|\, 3),\\
(1\,|\,2) & \mapsto (1\,|\,2).
\end{align*}
\item[(b)] Now let us take again $\xi=-1$ and write $Q=\xi^b$, $q=\xi^a$
where $a=1$ and $b=4$. Then we have $\cB=\Lambda_3^2$ and $\beta$
is the identity.
\end{itemize}

We leave it to the reader to extract this information from the
decomposition matrix printed in Table~\ref{tabb3}.
\end{exmp}

The following result was already announced in \cite[6.9]{mylaus}, with a
brief sketch of the proof. We include a more rigorous argument here.

\begin{thm} \label{asymp} Assume that $b>(n-1)a>0$. Then $\cB=
\Lambda_n^2$ is a canonical basic set, where the map $\beta
\colon \Lambda_n^2\rightarrow \cB$ is the identity.
\end{thm}

\begin{proof} First of all, the formula for $\ba_\ulambda$ in
Proposition~\ref{afuncWn} can be simplified under the given assumptions
on $a$ and $b$. Indeed, let $\ulambda \in \Pi_n^2$. Then, by
\cite[Example~3.6]{GeIa05}, we have
\begin{equation*}
f_{\ulambda}=1 \quad \mbox{and}\quad \ba_{\ulambda}=
b\,|\lambda^{(2)}|+a\,(n(\lambda^{(1)})+2n(\lambda^{(2)})-
n({\lambda^{(2)}}^*)),\tag{$*$}
\end{equation*}
where we set $n(\nu)=\sum_{i=1}^t (i-1)\nu_i$ for any partition $\nu=
(\nu_1\geq \nu_2 \geq \cdots \geq \nu_t\geq 0)$ and where $\nu^*$ denotes
the conjugate partition.  Thus, any field $k$ is $L$-good. The hypotheses
in Theorem~\ref{gero} concerning Lusztig's conjectures are satisfied by
\cite[Theorem~1.3]{GeIa05} and \cite[Corollary~7.12]{my05}. Thus, we already
know that a canonical basic set $\cB$ exists.

In order to prove that $\cB=\Lambda_n^2$, we must show that we have
the following implication, for any $\ulambda\in \Pi_n^2$ and $\umu\in
\Lambda_n^2$:
\[ [S^\ulambda:D^\umu] \neq 0 \quad \Rightarrow \quad
\ulambda=\umu \mbox{ or } \ba_{\umu}<\ba_{\ulambda}.\]
Using the relation ($\Delta$), we see that it is enough to prove
the following implication, for any $\ulambda,\mu\in \Pi_n^2$:
\[ \ulambda \trianglelefteq \umu \quad \Rightarrow \quad  \ulambda=\umu
\mbox{ or } \ba_\umu<\ba_\ulambda,\]
where $\trianglelefteq$ denotes the dominance order on bipartitions.
There are two cases.

{\em Case 1.} We have $|\lambda^{(1)}|=|\mu^{(1)}|$ and $|\lambda^{(2)}|=
|\mu^{(2)}|$. Then the condition $\ulambda \trianglelefteq \umu$ implies that
\[\lambda^{(1)} \trianglelefteq \mu^{(1)} \qquad \mbox{and}\qquad
\lambda^{(2)} \trianglelefteq \mu^{(2)},\]
where, on the right, the symbol $\trianglelefteq$ denotes the dominance
order on partitions. Now one can argue as in the proof of ``Case 1'' in
\cite[Corollary~5.5]{GeIa05} and conclude that $\ba_\umu\leq \ba_\ulambda$,
with equality only if $\umu= \ulambda$.

{\em Case 2.} We have $|\lambda^{(1)}|<|\mu^{(1)}|$ and $|\lambda^{(2)}|>
|\mu^{(2)}|$. If $b$ is very large with respect to $a$, then the formula
($*$) immediately shows that $\ba_\umu<\ba_\ulambda$. To get this conclusion
under the weaker condition that $b>(n-1)a$, one reduces to the case where
\begin{itemize}
\item $\lambda^{(1)}$ is obtained from $\mu^{(1)}$ by decreasing one
part by $1$, and
\item $\lambda^{(2)}$ is obtained from $\mu^{(2)}$ by increasing one
part by $1$.
\end{itemize}
(This is done by an argument similar to that in \cite[(I.1.16)]{Mac95}.)
Now one can argue as in the proof of ``Case 2'' in
\cite[Corollary~5.5]{GeIa05} and compare directly the values of $\ba_\umu$
and $\ba_\ulambda$. This yields $\ba_\umu<\ba_\ulambda$, as required.
\end{proof}

\section{The cases (A) and (B) in Theorem~1.1} \label{secAB}

In this section, we deal with cases (A) and (B) in
Theorem~\ref{mainspecht}, and show the existence of canonical basis sets.
Throughout, we fix a weight function $L\colon W_n \rightarrow \N$ such
that
\[L(t)=b \geq 0 \quad \mbox{and}\quad L(s_1)=\cdots =L(s_{n-1})=a\geq 0.\]
We also assume that $k$ is $L$-good; see Definition~\ref{goodL}. As before,
the parameters of $H_n$ are given by $Q=\xi^b$ and $q=\xi^a$, where $\xi
\in k^\times$. We will obtain the existence of canonical basic sets even
for those values of $a$ and $b$ where the hypotheses of Theorem~\ref{gero}
concerning Lusztig's conjectures are not yet known to hold.

\begin{thm} \label{caseA} Assume that $k$ is $L$-good and that we are in case
(A) of Theorem~\ref{mainspecht}, that is, we have $f_n(Q,q)\neq 0$. Then
\[\cB=\Lambda_n^2=\{\ulambda\in \Pi_n^2\mid \mbox{$\lambda^{(1)}$
and $\lambda^{(2)}$ are $e$-regular}\}\]
is a canonical basic set, where the map $\beta\colon \Lambda_n^2
\rightarrow \cB$ is the identity.
\end{thm}

\begin{proof} As already pointed out in Theorem~\ref{mainspecht}, we have
\[\Lambda_n^2=\{\ulambda\in \Pi_n^2\mid \mbox{$\lambda^{(1)}$
and $\lambda^{(2)}$ are $e$-regular}\}\]
under the given assumption on $f_n(Q,q)$.

Let us first deal with the case where $a=0$. Since $k$ is assumed to be
$L$-good, we have either $\mbox{char}(k)=0$ or $\mbox{char}(k)=p>n$.
In both cases, $e>n$. But then $H_n$ is semisimple and $D$ is the
identity matrix; see Dipper--James \cite[Theorem~5.5]{DJ2}. In particular,
$\cB=\Lambda_n^2=\Pi_n^2$ is a canonical basic set.

Let us now assume that $a>0$. We shall follow the argument in the proof
of \cite[Prop.~6.8]{my00}. By Dipper--James \cite{DJ2}, we can
express the decomposition numbers $[S^\ulambda:D^\umu]$ of
$H_n$ in terms of decomposition numbers for Iwahori--Hecke algebras
associated with the symmetric groups $\Sym_r$ for $0 \leq r \leq n$.

We need to set-up some notation. Let $H_k(\Sym_r,q)$ be the Iwahori--Hecke
algebra of $\Sym_r$, over the field $k$ and with parameter~$q$. By the
classical results of Dipper and James \cite{DJ0}, we have a Specht module
of $H_k(\Sym_r,q)$ for every partition $\lambda$ of $r$; let us denote this
Specht module by $S^\lambda$. (Again, our notation is such that $S^\lambda$
is the Specht module labelled by $\lambda^*$ in \cite{DJ0}.)
Furthermore, the simple modules of $H_k(\Sym_r, q)$ are labelled by the
$e$-regular partitions of $r$; let us denote by $D^\mu$ the simple module
labelled by the $e$-regular partition $\mu$ of~$r$. Correspondingly, we
have a matrix of composition multiplicities $[S^\lambda:D^\mu]$, such that
the following holds:
\begin{equation*}
\left\{\begin{array}{l} \quad [S^\mu:D^\mu]=1 \quad \mbox{for any
$e$-regular partition $\mu$ of $r$},\\ \quad [S^\lambda:D^\mu] \neq 0 \quad
\Rightarrow \quad \lambda \trianglelefteq \mu,\end{array}\right. \tag{$*$}
\end{equation*}
where $\trianglelefteq$ denotes the dominance order on partitions;
see \cite[Theorem~7.6]{DJ0}.

With this notation, we have the following result for $H_n$; see
Dipper--James \cite[Theorem~5.8]{DJ2}. Let $\ulambda= (\lambda^{(1)}\,|\,
\lambda^{(2)}) \in \Pi_n^2$ and $\umu=(\mu^{(1)}\,|\,\mu^{(2)})\in
\Lambda_n^2$. If $|\lambda^{(1)}|= |\mu^{(1)}|$ and
$|\lambda^{(2)}|= |\mu^{(2)}|$, then
\[ [S^\ulambda:D^\umu]=[S^{\lambda^{(1)}}:D^{\mu^{(1)}}]\cdot
[S^{\lambda^{(2)}}: D^{\mu^{(2)}}];\]
otherwise, we have $[S^\ulambda:D^\umu]=0$.

Now, in order to prove that $\Lambda_n^2$ is a canonical basic
set, we must show that we have the following implication:
\[ [S^\ulambda:D^\umu] \neq 0 \quad \Rightarrow \quad
 \ulambda=\umu \mbox{ or } \ba_{\umu}<\ba_{\ulambda}.\]
Taking into account the above results of Dipper and James, it will
be enough to consider bipartitions $\ulambda= (\lambda^{(1)}\,|\,
\lambda^{(2)})$ and $\umu=(\mu^{(1)}\,|\,\mu^{(2)})$ such that
$|\lambda^{(1)}|=|\mu^{(1)}|$ and $|\lambda^{(2)}|= |\mu^{(2)}|$. For
such bipartitions, we must show the following implication:
\[ [S^{\lambda^{(1)}}:D^{\mu^{(1)}}]\neq 0 \mbox{ and }
[S^{\lambda^{(2)}}:D^{\mu^{(2)}}]\neq 0 \quad \Rightarrow \quad
\ulambda=\umu \mbox{ or } \ba_{\umu}<\ba_{\ulambda}.\]
Using ($*$), it will actually be sufficient to prove the following
implication:
\begin{equation*}
\lambda^{(1)} \domi \mu^{(1)} \mbox{ and }
 \lambda^{(2)} \domi \mu^{(2)} \quad \Rightarrow \quad
\ulambda=\umu \mbox{ or } \ba_{\umu}<\ba_{\ulambda}.\tag{$\dagger$}
\end{equation*}
Thus, we are reduced to a purely combinatorial statement. We claim that,
in fact, ($\dagger$) holds for all bipartitions $\ulambda,\umu\in
\Pi_n^2$. To prove this, one can further reduce to the situation where
$\lambda^{(1)}=\mu^{(1)}$ or $\lambda^{(2)}=\mu^{(2)}$. For example, let
us assume that $\lambda^{(1)}=\mu^{(1)}$ and $\lambda^{(2)} \domi
\mu^{(2)}$. Then one can even further reduce to the case where $\lambda^{(2)}$
is obtained from $\mu^{(2)}$ by increasing one part by $1$ and by
decreasing another part by~$1$; see Macdonald \cite[(I.1.16)]{Mac95}. But
in this situation, it is straightforward to check the desired assertion
by directly using the formula in Proposition~\ref{afuncWn}. Thus,
($\dagger$) is proved, and this completes the proof that
$\cB=\Lambda_n^2$ is a canonical basic set.
\end{proof}

Now let us turn to case (B) in Theorem~\ref{mainspecht}. Thus, we assume
that $f_n(Q,q)=0$ and $q=1$; note that this also gives $Q=-1$. Furthermore,
we have
\[ e=\left\{\begin{array}{cl} \infty & \qquad \mbox{if $\mbox{char}(k)=0$},
\\ p & \qquad \mbox{if $\mbox{char}(k)=p>0$}.\end{array}\right.\]
Note also that, since $k$ is assumed to be $L$-good, we must have
$b>0$ and $a$ cannot divide~$b$.

We need some results concerning the complex irreducible characters of
the group $W_n$. These characters are naturally labelled by $\Pi_n^2$;
we write $\chi^\ulambda$ for the irreducible character labelled by
$\ulambda \in \Pi_n^2$. (In fact, as already noted earlier, the simple
modules of $\bH_{K,n}$ are given by the various Specht modules $S^\ulambda$;
then $\chi^\ulambda$ is the character of the corresponding Specht
module of $W_n$ obtained by specializing $u\mapsto 1$.) Now consider the
parabolic subgroup
\[\Sym_n=\langle s_1,\ldots,s_{n-1}\rangle \subseteq W_n.\]
Its irreducible characters are labelled by the set $\Pi_n$ of all
partitions of $n$; let us denote by $\psi^\nu$ the irreducible character
labelled by $\nu\in \Pi_n$. We set
\[ \ba_\nu=\sum_{i=1}^t a(i-1)\nu_i \qquad
\mbox{if $\nu=(\nu_1\geq \nu_2\geq  \cdots \geq \nu_t\geq 0)$}.\]
(This is the $\ba$-function for the Iwahori--Hecke algebra of the
symmetric group with weight function given by the constant value $a$ on
each generator; see \cite[22.4]{Lusztig03}.) Now, by general properties
of Lusztig's $\ba$-function, there is a well-defined map (called
``truncated induction'' or ``$J$-induction'')
\[ J\colon \Pi_{n} \rightarrow \Pi_n^2\]
such that $\ba_{J(\nu)}=\ba_\nu$ and
\[ \Ind_{\Sym_n}^{W_n}(\psi^\nu)=\chi^{J(\nu)}+\mbox{sum of
characters $\chi^{\ulambda}$ such that $\ba_{\ulambda}>\ba_{J(\nu)}$}.\]
(This could be deduced from the results in \cite[Chap.~20]{Lusztig03}, but
Lusztig assumes the validity of (P1)--(P15) in that chapter. A reference
which does not refer to these properties is provided by \cite[\S 3]{my02}.)

\begin{rem} \label{aa0} Assume that $a=0$. Then $\ba_\nu=0$ for all $\nu
\in \Pi_n$. A formula for $\ba_\ulambda$ (where $\ulambda \in \Pi_n^2$)
is given in Remark~\ref{expa}. One easily checks that
\[ J(\nu)=(\nu,\varnothing) \qquad \mbox{for all $\nu \in \Pi_n$}.\]
Now assume that $a>0$ and write $b=ar+b'$ as before; we have already
noted above that $b'>0$ (since $k$ is assumed to be $L$-good).
Then one can give a combinatorial description of the map $J \colon \Pi_n
\rightarrow \Pi_n^2$, using the results in \cite[Chapter~22]{Lusztig03}.
Indeed, first note that $\chi^\nu$ is obtained by $J$-inducing the sign
character of the Young subgroup $\Sym_{\nu^*}$ to $\Sym_n$, where $\nu^*$
denotes the dual partition. Thus, we have
\[ \Ind_{\Sym_{\nu^*}}^{W_n}(\mbox{sign})=\chi^{J(\nu)}+\mbox{sum of
characters $\chi^{\ulambda}$ such that $\ba_{\ulambda}>\ba_{J(\nu)}$}.\]
Now write $\nu^*=(\nu_1^*>\nu_2^*>\ldots>\nu_t^*>0)$. Then the induction
can be done in a sequence of steps, using the following embeddings:
\begin{align*}
\Sym_{\nu_1^*} &\subseteq W_{\nu_1^*}, \\ W_{\nu_1^*}\times
\Sym_{\nu_2^*} & \subseteq W_{\nu_1^*+\nu_2^*},\\ W_{\nu_1^*+\nu_2^*}
\times \Sym_{\nu_3^*} & \subseteq W_{\nu_1^*+\nu_2^* +\nu_3^*},\\ \vdots &
\\ W_{\nu_1^*+\nu_2^*+\cdots +\nu_{t-1}^*}\times \Sym_{\nu_t^*} &
\subseteq W_{n}.
\end{align*}
The combinatorial rule for describing the $J$-induction in the above
steps is provided by \cite[Lemma~22.17]{Lusztig03}. For this purpose,
choose a large integer $N$ and consider the symbol
\[ X_0=\left(\begin{array}{cc} a+b',\; 2a+b',\;3a+b',\;\ldots,\;(N+r)a+b'\\
a,\quad 2a,\quad 3a,\quad \ldots,\quad Na\end{array}\right).\]
(This labels the trivial character of the group $W_0=\{1\}$.) We define a
sequence of symbols $X_1,X_2,\ldots,X_t$ as follows. Let $X_1$ be the symbol
obtained from $X_0$ by increasing each of the $\nu_1^*$ largest entries in
$X_0$ by~$a$. Then let $X_2$ be the symbol obtained from $X_1$ by increasing
each of the $\nu_2^*$ largest entries in $X_1$  by~$a$, and so on. Then the
symbol $X_t$ corresponds, by the reverse of the procedure described in
Section~2, to the bipartition $J(\nu)$.
\end{rem}

\begin{exmp} \label{trunc} (i) Assume that $a>0$ and $r\geq n-1$. Then,
in each step of the construction of $X_t$, we will only increase entries in
the top row of a symbol. Thus, we obtain
\[ J(\nu)=(\nu,\varnothing) \qquad \mbox{for all $\nu \in \Pi_n$}.\]
(ii) In general, let $\nu\in \Pi_n$ and suppose that $J(\nu)=(\lambda\,|\,
\mu)$ where $\lambda$ and $\mu$ are partitions. Then one can check that
$\nu$ is the union of all the parts of $\lambda$ and $\mu$ (reordered if
necessary). For example, assume that $a=2$ and $b=1$; then $r=0$ and $b'=1$.
We obtain
\begin{align*}
J(4)&=(4\,|\,\varnothing),\\ J(31)&=(3\,|\,1),\\ J(22)&=(2\,|\,2),\\
J(211)&=(21\,|\,1),\\ J(1111)&=(11\,|\,11).
\end{align*}
\end{exmp}

\begin{thm} \label{caseB} Assume that $k$ is $L$-good and that we are
in case (B) of Theorem~\ref{mainspecht}, where
\[ \Lambda_n^2=\{\ulambda \in \Pi_n^2\mid \mbox{$\lambda^{(1)}$
is $e$-regular and $\lambda^{(2)}=\varnothing$}\}.\]
Then $\cB=\beta(\Lambda_n^2)$ is a canonical basic set for $H_n$,
where $\beta\colon \Lambda_n^2 \rightarrow \Pi_n^2$ is defined
by $(\lambda^{(1)}\,|\,\varnothing) \mapsto J(\lambda^{(1)})$.
\end{thm}

\begin{proof} The special feature of case (B) is that, by Dipper--James
\cite[5.4]{DJ2}, the simple $H_n$-modules are obtained by extending (in a
unique way) the simple modules of the parabolic subalgebra $k\Sym_n=
\langle T_{s_1}, \ldots,T_{s_{n-1}} \rangle_k$ to $H_n$. Thus, given
an $e$-regular partition $\nu\in \Pi_n$, the restriction
\[ D^\nu:=\mbox{Res}^{H_n}_{k\Sym_n}(D^{(\nu\,|\,\varnothing)})\]
is a simple module for $k\Sym_n$, and all simple $k\Sym_n$-modules
are obtained in this way. Moreover, this yields precisely the classical
labelling of the simple $k\Sym_n$-modules by $e$-regular partitions of $n$.

Alternatively, we can express this by saying that the projective
indecomposable $H_n$-modules are obtained by inducing the projective
indecomposable $k\Sym_n$-modules to $H_n$. More precisely, for any $e$-regular
$\nu \in \Pi_n$, let $P^\nu$ be a projective cover of $D^\nu$; then
\[ Q^\nu:=\Ind_{k\Sym_n}^{H_n}(P^\nu) \]
is a projective cover of the simple $H_n$-module $D^{(\nu \,
|\,\varnothing)}$. (This follows from standard results on projective
modules and Frobenius reciprocity.)

It will now be convenient to use an appropriate form of Brauer recipricity
in the usual setting of modular representation theory (as in
\cite[Remark~5.9]{DJ2}). Then the multiplicity of a simple module in a
Specht module is seen to be the same as the multiplicity (in the appropriate
Grothendieck group) of that Specht module in the projective cover of the
simple module. Thus, for any $e$-regular $\nu \in \Pi_n$, we have
\[ [P^\nu]=\sum_{\lambda\in \Pi_n} [S^\lambda:D^\nu]\cdot [S^\lambda],\]
where $[P^\nu]$, $[S^\lambda]$ denote the classes of these modules in the
appropriate Grothendieck group of $k\Sym_n$-modules, and we have
\[ [Q^\nu]=\sum_{\ulambda\in \Pi_n^2} [S^\ulambda:D^{(\nu\,|\,
\varnothing)}] \cdot [S^{\ulambda}],\]
where $[Q^\nu]$, $[S^\ulambda]$ denote the classes of these modules in
the appropriate Grothendieck group of $H_n$-modules.

With these preparations, let us now consider the $\ba$-invariants.
Using the relations ($*$) in the proof of Theorem~\ref{caseA}, we have:
\[ [P^\nu]=[S^\nu]+\mbox{lower terms with respect to $\domi$},\]
where ``lower terms'' stands for a sum of classes of modules $S^{\nu'}$
such that $\nu' \domi \nu$ and $\nu'\neq \nu$. By known properties of the
$\ba$-function, this can also be expressed as:
\[ [P^\nu]=[S^\nu]+\mbox{higher terms},\]
where ``higher terms'' stands for a sum of classes of modules $S^{\nu'}$
such that $\ba_{\nu'}> \ba_\nu$. Now let us induce to $H_n$. By the
definition of $J(\nu)$, we have
\[[\Ind_{k\Sym_n}^{H_n}(S^\nu)]=[S^{J(\nu)}]+\mbox{higher terms},\]
where $\ba_\nu=\ba_{J(\nu)}$ and where ``higher terms'' stands for a sum
of classes of modules $S^{\nu'}$ such that $\ba_{\nu'}> \ba_\nu$.
Furthermore, by general properties of the $\ba$-function, it is known that
inducing a module with a given $\ba$-invariant $i$ will result in a sum of
modules of $\ba$-invariants $\geq i$. Hence we obtain
\begin{align*}
[Q^\nu]&=[\Ind_{k\Sym_n}^{H_n}(S^\nu)]
+\mbox{higher terms}\\ &=[S^{J(\nu)}]+\mbox{higher terms},
\end{align*}
where, in both cases,  ``higher terms'' stands for a sum of classes of
modules $S^{\unu'}$ such that $\ba_{\unu'}>\ba_\nu=\ba_{J(\nu)}$.
Thus, the conditions ($\Delta_\ba$) in Definition~\ref{defbs} are satisfied
with respect to the map
\[ \beta\colon \Lambda_n^2 \rightarrow \Pi_n^2, \qquad (\nu\,|\,\varnothing)
\mapsto J(\nu).\]
Hence, $\cB=\{J(\nu)\mid \mbox{$\nu\in \Pi_n$ is $e$-regular}\}$
is a canonical basic set.
\end{proof}

\section{The Fock space and canonical bases} \label{sec:klesh}

In this section, we briefly review the deep results of Ariki and Uglov
concerning the connections between the representation theory of Hecke
algebras and the theory of canonical bases for quantum groups. (For general
introductions to the theory of canonical bases, see Kashiwara \cite{Kash}
and Lusztig \cite{Lucan}.) These results will be used in the subsequent
section to describe the canonical basic set in case (C) of
Theorem~\ref{mainspecht}. The main theorems of this section are available
for a wider class of algebras, namely, the Arike--Koike algebras which
we now define.

\begin{defn} \label{arko} Let $k$ be an algebraically closed field and let
$\zeta\in k^\times$. Let $n,r\geq 1$ and fix parameters
\[ \bu=(u_1,\ldots,u_r) \qquad \mbox{where} \qquad  u_i\in \Z.\]
Having fixed these data, we let $H_{n,\zeta}^{\bu}$ be the associative
$k$-algebra (with $1$), with generators $S_0,S_1,\ldots, S_{n-1}$ and defining
relations as follows:
\begin{gather*} S_0S_1S_0S_1=S_1S_0S_1S_0 \quad \mbox{and}\quad
S_0S_i=S_iS_0 \quad \mbox{(for $i>1$)},\\ S_iS_j=S_jS_i \quad \mbox{(if
$|i-j|>1$)}, \\ S_iS_{i+1} S_i=S_{i+1}S_iS_{i+1} \quad \mbox{(for $1 \leq i
\leq n-2$)},\\ (S_0-\zeta^{u_1})(S_0- \zeta^{u_2})\cdots
(S_0-\zeta^{u_r})=0,\\ (S_i-\zeta)(S_i+1)=0 \quad \mbox{for $1 \leq i \leq
n-1$)}.
\end{gather*}
This algebra can be seen as an Iwahori--Hecke algebra associated with the
complex reflection group $G_{r,n}:=(\Z/r\Z)^n \rtimes \Sym_n$. See
Ariki \cite[Chap.~13]{Ar3} and Brou\'e--Malle \cite{BM2} for further
details and motivations for studying this class of algebras.
\end{defn}

\begin{rem} \label{ident} Let $r=2$. Then we can identify
$H_{n,\zeta}^{\bu}$ with an Iwahori--Hecke algebra of type $B_n$. Indeed,
the generator $S_0$ satisfies the quadratic relation
\[ (S_0-\zeta^{u_1})(S_0-\zeta^{u_2})=0.\]
Then the map $T_t \mapsto -\zeta^{-u_2}S_0$, $T_{s_1}\mapsto S_1$, $\ldots$,
$T_{s_{n-1}}\mapsto S_{n-1}$ defines an isomorphism
\[ H_k(W_n,-\zeta^{u_1-u_2},\zeta) \stackrel{\sim}{\rightarrow}
H_{2,\zeta}^{\bu}.\]
Note that, if $\zeta \neq 1$ and if $u_1,u_2$ are such that $|u_1-u_2|
\leq n-1$, then $f_n(-\zeta^{u_1-u_2},\zeta)=0$ and we are in Case~(C) of
Theorem~\ref{mainspecht}.
\end{rem}

Now Dipper--James--Mathas \cite[\S 3]{DJMa} have generalized the theory of
 Specht modules to the algebras $H_{n,\zeta}^{\bu}$; see also Graham--Lehrer
\cite{GrLe}. Let $\Pi_n^r$ denote the set of all ($r$-)multipartitions
of $n$, that is, $r$-tuples of partitions $\ulambda=(\lambda^{(1)}\,
|\,\ldots\,|\, \lambda^{(r)})$ such that $|\lambda^{(1)}|+\cdots +
|\lambda^{(r)}|=n$. For any $\ulambda\in \Pi_n^r$, there is a {\em Specht
module} $S^{\ulambda,\bu}$ for $H_{n,\zeta}^{\bu}$. Each $S^{\ulambda,\bu}$
carries a symmetric bilinear form and, taking quotients by the radical, we
obtain a collection of modules $D^{\ulambda, \bu}$. As before, we set
\[ \Lambda_n^{\bu}:=\{\ulambda \in \Pi_n^r\mid D^{\ulambda,\bu}
\neq \{0\}\}.\]
Then, by \cite[Theorem~3.30]{DJMa}, we have
\[ \Irr(H_{n,\zeta}^{\bu})=\{D^{\ulambda,\bu} \mid \ulambda \in
\Lambda_{n}^{\bu}\}.\]
Furthermore, the entries of the decomposition matrix
\[ D=\bigl([S^{\ulambda,\bu} :D^{\umu,\bu}]\bigr)_{\ulambda \in \Pi_n^r,
\umu \in \Lambda_{n}^{\bu}}\]
satisfy the conditions
\begin{equation*}
\left\{\begin{array}{l} \quad [S^{\umu,\bu}:D^{\umu,\bu}]=1 \quad
\mbox{for any $\umu\in \Lambda_{n}^{\bu}$},\\ \quad
[S^{\ulambda,\bu}:D^{\umu,\bu}] \neq 0 \quad \Rightarrow \quad
\ulambda \trianglelefteq \umu,\end{array}\right.  \tag{$\Delta^{\bu}$}
\end{equation*}
where $\trianglelefteq$ denotes the dominance order on $r$-partitions,
as defined in \cite[3.11]{DJMa}. Note, again, that these conditions
uniquely determine the set $\Lambda_{n}^{\bu}$ once the matrix $D$ is known.
By Ariki \cite{Ar2}, the problem of computing $D$ (at least in the case
where $\mbox{char}(k)=0$) can be translated to that of computing the
canonical bases of a certain module over the quantum group
$\cU_q(\widehat{\fsl}_e)$, where $e\geq 2$ is the order of $\zeta$ in
the multiplicative group of~$k$.

We first give a brief overview of the results of Uglov \cite{Ug} which
generalize previous work of Leclerc and Thibon \cite{LT}; for a good survey
on this theory, see Yvonne \cite{Yv}. Let us fix an integer $e\geq 2$. Let
$\bu=(u_1,\ldots,u_r)\in\Z^r$ and let $q$ be an indeterminate.
The Fock space $\cF^{\bu}$ is defined to be the $\C(q)$-vector space
generated by the symbols $|\ulambda,\bu\rangle$ with $\ulambda\in \Pi_n^r$:
\[ \cF^{\bu}:=\bigoplus_{n=0}^{\infty} \bigoplus_{\ulambda \in\Pi_n^r}
\C(q)\,|\,\ulambda,\bu\rangle\]
where $\Pi_0^r=\{\uvar=(\varnothing,\ldots,\varnothing)\}$.
Let $\cU_q^\prime(\widehat{\fsl}_e)$ be the quantum group associated to
the Lie algebra $\widehat{\fsl}_e^\prime$. Then the deep results of Uglov
show how the set $\cF^{\bu}$ can be endowed with a structure of integrable
$\cU_q^\prime({\widehat{\fsl}}_e)$ module (see \cite[\S 3.5, \S 4.2]{Ug}).
Moreover, Uglov has defined an involution $\bar{\;}\;\colon \cF^{\bu}
\rightarrow \cF^{\bu}$. Then one can show that there is a unique basis
\[ \{G(\ulambda,\bu)\ |\ \ulambda\in\Pi_n^r,\ n\in \N\}\]
of $\cF^{\bu}$ such that the following two conditions hold:
\begin{align*}
\overline{ G(\ulambda,\bu)}& =G(\ulambda,\bu),\\
G(\ulambda,\bu)&=|\ulambda,\bu\rangle+ \mbox{$q\C[q]$-combination of
basis elements  $|\umu,\bu\rangle$}.
\end{align*}
The set $\{G(\ulambda,\bu)\}$ is called the Kashiwara--Lusztig
{\em canonical basis} of $\cF^{\bu}$.

Now we consider the $\cU_q^\prime(\widehat{\fsl_e})$-submodule $\cM^{\bu}
\subseteq \cF^{\bu}$ generated by $|\uvar,\bu\rangle$. It is well-known
that this is isomorphic to the irreducible
$\cU_q^\prime(\widehat{\fsl}_e)$-module $V(\Lambda)$ with highest weight
\[ \Lambda:=\Lambda_{u_1(\bmod e)}+\Lambda_{u_2(\bmod e)}+\ldots+
\Lambda_{u_r(\bmod e)}.\]
A basis of $\cM^{\bu}$ can be given by using the canonical
basis of $\cF^{\bu}$ and by studying the associated crystal graph. To
describe this graph, we will need some further combinatorial definitions.

Let $\ulambda= (\lambda^{(1)}\,|\, \ldots\,|\,\lambda^{(r)})\in \Pi_n^r$
and write
\[ \lambda^{(c)}=(\lambda^{(c)}_1 \geq \lambda^{(c)}_2 \geq \cdots
\geq 0) \qquad \mbox{for} \qquad c=1,\ldots,r.\]
The diagram of $\ulambda$ is defined as the set
\[ [\ulambda]:=\{(a,b,c) \mid 1 \leq c \leq r, \;1 \leq b
\leq \lambda^{(c)}_a \mbox{ for $a=1,2,\ldots$}\}.\]
For any ``node'' $\gamma=(a,b,c)\in [\ulambda]$, we set
\[ \mbox{res}_e(\gamma):=(b-a+u_c)  \quad \bmod e\]
and call this the {\em $e$-residue} of $\gamma$ with respect to the
parameters $\bu$. If $\mbox{res}_e(\gamma)=i$, we say that
$\gamma$ is an $i$-node of $\ulambda$.

Now suppose that $\ulambda \in \Pi_n^r$ and $\umu\in \Pi_{n+1}^r$ for
some $n\geq 0$. We write
\[ \gamma=\umu/\ulambda \qquad \mbox{if}\qquad
[\ulambda] \subset [\umu] \quad \mbox{and}\quad
[\umu]=[\ulambda] \cup \{\gamma\};\]
Then we call $\gamma$ an addable node for $\ulambda$ or a
removable node for $\umu$.

\begin{defn}[Foda et al. \protect{\cite[p.~331]{FLOTW}}] \label{above}
We say that the node $\gamma=(a,b,c)$ is ``above'' the node $\gamma'=(a',
b',c')$ if
\begin{itemize}
\item either $b-a+u_c<b'-a'+u_{c'}$,
\item or $b-a+u_c=b'-a'+u_{c'}$ and $c'<c$.
\end{itemize}
\end{defn}

Using this order relation on nodes, we define the notion of ``good'' nodes,
as follows. Let $\ulambda\in \Pi_n^r$ and let $\gamma$ be an $i$-node of
$\ulambda$. We say that $\gamma$ is a {\em normal} node if, whenever
$\gamma'$ is an $i$-node of $\ulambda$ below $\gamma$, there are strictly
more removable $i$-nodes between $\gamma'$ and $\gamma$ than there are
addable $i$-nodes between $\gamma'$ and $\gamma$. If $\gamma$ is a highest
normal $i$-node of $\ulambda$, then $\gamma$ is called a {\em good} node.
Note that these notions heavily depend on the definition of what it means
for one node to be ``above'' another node. These definitions (for $r=1$)
first appeared in the work of Kleshchev \cite{Klesh} on the modular
branching rule for the symmetric group; see also the discussion of these
results in \cite[\S 2]{LLT}.

\begin{defn} \label{defuglov} For any $n\geq 0$, we define a subset
$\Phi_{e,n}^{\bu} \subseteq \Pi_n^r$ recursively as follows. We set
$\Phi_{e,0}^{\bu}=\{\uvar\}$. For $n\geq 1$, the set $\Phi_{e,n}^{\bu}$ is
constructed as follows.
\begin{itemize}
\item[(1)] We have $\uvar \in \Phi_{e,n}^{\bu}$;
\item[(2)] Let $\ulambda\in \Pi_n^r$. Then $\ulambda$
belongs to $\Phi_{e,n}^{\bu}$ if and only if $\ulambda/\umu=\gamma$ where
$\umu \in \Phi_{e,n-1}^{\bu}$ and $\gamma$ is a good $i$-node of $\ulambda$
for some $i\in\{0,1,\ldots,e-1\}$.
\end{itemize}
The set
\[ \Phi_e^{\bu}:=\bigcup_{n\geq 0} \Phi_{e,n}^{\bu}\]
will be called the set of {\em Uglov $r$-multipartitions}.
\end{defn}

\begin{rem} \label{invar} Let $m \in Z$ and consider the parameter set
$\bu'=(u_1+m,u_2+m,\ldots,u_r+m)$. Note  that a node $\gamma$ is
above a node $\gamma'$ with respect to $\bu$ if and only if this holds
with respect to $\bu'$. It follows that $\Phi_e^{\bu}=\Phi_e^{\bu'}$.
\end{rem}

\begin{thm}[Jimbo et al. \protect{\cite{JMMO}}, Foda et al.
\protect{\cite{FLOTW}}, Uglov \protect{\cite{Ug}}] \label{jmmo} The crystal
graph of $\cM^{\bu}$ has vertices labelled by the  set $\Phi_e^{\bu}$ of Uglov
$r$-multipartitions. Given two vertices $\ulambda\neq \umu$ in that graph, we
have an edge
\[ {\ulambda} \stackrel{i}{\rightarrow} {\umu} \qquad \mbox{(where
$0 \leq i \leq e-1$)}\]
if and only if $\umu$ is obtained from $\ulambda$ by adding a  ``good''
$i$-node
\end{thm}

For a general introduction to crystal graphs, see Kashiwara \cite{Kash}.

\begin{rem} \label{kleshdef} (a) Assume that  $\bu\in \Z^r$ is such that
\[ 0\leq u_1 \leq u_2\leq  \ldots \leq u_{r}\leq e-1.\]
Then it is shown in Foda et al. \cite[2.11]{FLOTW} that $\ulambda \in
\Pi_{r,n}$ belongs to $\Phi_{e,n}^{\bu}$ if and only if the following
conditions are satisfied:
\begin{itemize}
\item For all $1\leq j \leq r-1$ and $i=1,2,\ldots$, we have:
\[ \lambda^{(j+1)}_i \geq \lambda^{(j)}_{i+u_{j+1}-u_j} \qquad \mbox{and}
\qquad \lambda^{(1)}_i \geq \lambda^{(r)}_{i+l+u_1-u_r};\]
\item for all $k>0$, among the residues appearing at the right ends
of the rows of $[\ulambda]$ of length $k$, at least one element of
$\{0,1,\ldots,e-1\}$ does not occur.
\end{itemize}
Note that this provides a non-recursive description of the elements
of $\Lambda_n^{\bu}$.

(b) Assume that $\bu\in \Z^r$ is such that
\[ u_1>u_2 > \cdots > u_r>0 \qquad \mbox{where}\qquad u_i-u_{i+1}>n-1
\mbox{ for all $i$}.\] Then the set of Uglov $r$-multipartitions
$\Phi_{e,n}^{\bu}$ coincides with the set $\cK_{e,n}^{\bu}$ of Kleshchev
$r$-multipartitions as defined by Ariki \cite{Ar2}. More directly,
$\cK_{e,n}^{\bu}$ can be defined recursively in a similar way as in
Definition~\ref{defuglov}, where we use the following order relation on
nodes $\gamma=(a,b,c)$ and $\gamma=(a',b',c')$:
\begin{equation*}
\gamma \mbox{ is above } \gamma' \qquad
\stackrel{\text{def}}{\Leftrightarrow} \qquad c'<c \quad \mbox{or if}
\quad c=c' \mbox{ and } a'<a.
\end{equation*}
This is the order on nodes used by Ariki \cite[Theorem~10.10]{Ar3}.
\end{rem}

Finally, the following theorem gives a link between the canonical basis
elements of $\cM^{\bu}$ and the decomposition matrices of Ariki-Koike
algebras. The previous results show that the canonical basis of
$\cM^{\bu}$ is given by
\[ \{G(\umu,\bu)\, |\, \umu\in\Phi_e^{\bu}\}.\]
For each $\umu\in\Phi_{e,n}^{\bu}$, we can write
\[ G(\umu,\bu)=\sum_{\ulambda\in \Pi_n^r} d_{\ulambda,\umu}^\bu(q)\,
|\ulambda,\bu\rangle \qquad \mbox{where $d_{\ulambda,\umu}^\bu(q)
\in\C[q]$}.\]
With this notation, we can now state:

\begin{thm}[Ariki \protect{\cite{Ar},\cite{Ar2},\cite{Ar3}}]
\label{thmariki} Let $e\geq 2$, $\bu\in\Z^r$ and let $H_{n,\zeta}^{\bu}$
be the Ariki-Koike algebra over $k$ as in Definition~\ref{arko}, where
$1 \neq \zeta \in k^\times$ is a root of unity and $k$ has characteristic
$0$. Let $e\geq 2$ be the order of $\zeta$. By adding multiples of $e$
to each $u_i$, we may assume without loss of generality that
\[ u_1>u_2 > \cdots > u_r>0 \qquad \mbox{where}\qquad u_i-u_{i+1}>n-1
\mbox{ for all $i$}.\]
Then $\Phi_{e,n}^{\bu}=\cK_{e,n}^{\bu}=\Lambda_{n}^\bu$ and
\[[S^{\ulambda,\bu}:D^{\umu,\bu}]=d_{\ulambda,\umu}^\bu(1)\]
for all $\ulambda\in \Pi_n^r$ and  $\umu\in\Phi_{e,n}^{\bu}$.
\end{thm}

As a consequence, we can compute the decomposition matrices of Ariki--Koike
algebras using the known combinatorial algorithm for computing canonical basis
for $\cM^{\bu}$; see Lascoux--Leclerc--Thibon \cite{LLT} for $r=1$, and
Jacon \cite{Jac3} for $r\geq 2$.

\begin{cor} \label{otherui} Let us keep the same general hypotheses as in
Theorem~\ref{thmariki}, except that we drop the condition that
$u_i-u_{i+1}>n-1$ for all~$i$. Then there exists a bijection $\kappa
\colon \Lambda_n^\bu\rightarrow \Phi_{e,n}^{\bu}$ such that
\[[S^{\ulambda,\bu}:D^{\umu,\bu}]=d_{\ulambda,\kappa(\umu)}^\bu(1)\]
for all $\ulambda\in \Pi_n^r$ and  $\umu\in\Lambda_n^\bu$.
\end{cor}

\begin{proof} Let $\bu'=(u_1',\ldots,u_r')\in \Z^r$ be such that
$u_i'-u_{i+1}'>n-1$ and $u_i\equiv u_i'\bmod e$ for all~$i$.
As explained in Foda et al. \cite[Note~2.7]{FLOTW}, the canonical basis
$\{G(\umu,\bu)\}$ (specialised at $q=1$) coincides with the canonical
basis $\{G(\umu,\bu')\}$ (specialised at $q=1$), at least as far as all
multipartitions $\umu$ of total size $\leq n$ are concerned. Hence there
exists a bijection $\kappa \colon \Phi_{e,n}^{\bu'} \rightarrow
\Phi_{e,n}^{\bu}$ such that
\[ d_{\ulambda,\kappa(\umu)}^{\bu}(1)=d_{\ulambda,\umu}^{\bu'}(1) \quad
\mbox{for all $\ulambda\in \Pi_n^r$ and  $\umu\in\Phi_{e,n}^{\bu}$}.\]
By Theorem~\ref{thmariki}, we have $\Phi_{e,n}^{\bu'}=
\Lambda_n^{\bu'}$ and $[S^{\ulambda,\bu'}:D^{\umu,\bu'}]=
d_{\ulambda,\umu}^{\bu'}(1)$. Now note that $H_{n,\zeta}^\bu=
H_{n,\zeta}^{\bu'}$ and so $\Lambda_n^{\bu}= \Lambda_n^{\bu'}$. Thus, we
obtain
\[ [S^{\ulambda,\bu}:D^{\umu,\bu}]=
[S^{\ulambda,\bu'}:D^{\umu,\bu'}]=d_{\ulambda,\umu}^{\bu'}
(1)=d_{\ulambda,\kappa(\umu)}^{\bu}(1)\]
for all $\ulambda\in \Pi_n^r$ and  $\umu\in\Lambda_{n}^{\bu}$, as required.
\end{proof}

\section{The case (C) in Theorem~1.1} \label{secC}

Using the results of the previous section, we are now ready to deal with
case (C) in Theorem~\ref{mainspecht} and show the existence of canonical
basis sets when $k$ has characteristic zero. (The case of positive
characteristic remains conjectural, see Remark~\ref{reduct1}.) Throughout,
we fix a weight function $L\colon W_n \rightarrow \N$ as in
Section~\ref{secAB}. As before, the parameters of $H_n$ are given by 
$Q=\xi^b$ and $q=\xi^a$, where $\xi \in k^\times$ and $a,b\geq 0$. We 
shall now assume that
\begin{itemize}
\item[(C1)] $Q=-q^d$ for some $d \in \Z$,
\item[(C2)] $q\neq 1$ has finite order in $k^\times$,
\item[(C3)] $\mbox{char}(k)\neq 2$.
\end{itemize}

\begin{rem} \label{C0} Assume that $\mbox{char}(k)\neq 2$, $f_n(Q,q)=0$ 
and $q\neq 1$. Then $\xi$ is a non-trivial element of finite even order in 
the multiplicative group of~$k$; furthermore, (C1), (C2) and (C3) hold. 

Indeed, the condition $f_n(Q,q)=0$ implies that $Q=-q^d$ for some
integer~$d$ such that $|d|\leq n-1$. Thus, we have
\[ \xi^{b-ad}=-1 \qquad \mbox{where $-(n-1)\leq d \leq n-1$}.\]
Now, if $b\neq ad$, then this relation shows that $\xi$ is a non-trivial
element of finite even order in $k^\times$. If we had $b=ad$, we would
obtain the contradiction $1=\xi^{b-ad}=-1$. Thus, the above claim is proved.
\end{rem}

Thus, assuming that $\mbox{char}(k)\neq 2$, the conditions (C1), (C2) are 
somewhat weaker than the condition in case~(C) of Theorem~1.1, as $d$ is 
not required to be of absolute value $\leq n-1$. The following results 
will hold assuming only (C1), (C2), (C3). We have seen above that $\xi$ has 
finite order; let $l\geq 2$ be the multiplicative order of $\xi$. Let 
$\zeta_l\in \C$ be a primitive $l$-th root of unity. We shall consider 
the Iwahori--Hecke algebra
\[ H_n^0:=H_{\C}(W_n,Q_0,q_0) \qquad \mbox{where} \qquad Q_0:=\zeta_l^b,
\qquad q_0:=\zeta_l^a \neq 1.\]
Note that both $H_n$ and $H_n^0$ are obtained by specialisation from
the same generic algebra $\bH_n$ (defined with respect to the given
weight function $L$). Note also that $f_n(Q,q)=0 \Leftrightarrow 
f_n(Q_0,q_0)=0$, and that the parameter $e$ defined with respect to 
$q_0$ is the same as the parameter $e$ defined with respect to~$q$ 
(since $q\neq 1$). Now Theorem~\ref{mainspecht} shows that the 
simple modules of $H_n$ and of $H_n^0$ are both parametrized by the same 
set $\Lambda_{n}^2$. In particular, we have $|\Irr(H_n)|=|\Irr(H_n^0)|$; 
see also Ariki--Mathas \cite[Theorem~A]{ArMa}. We have the following 
result, which reduces the determination of a canonical basic set to the 
case where $k$ has characteristic zero (assuming that such basic
sets exist at all).

\begin{lem}[See \protect{\cite[\S 3.1B]{Jac0}}] \label{reduct0} If $H_n$
admits a canonical basic set $\cB$ (with respect to a map $\beta\colon
\Lambda_n^2 \rightarrow \Pi_n^2$) and $H_n^0$ admits a canonical
basic set $\cB^0$ (with respect to a map $\beta^0\colon \Lambda_n^2
\rightarrow \Pi_n^2$), then we have $\cB=\cB^0$ and $\beta=\beta^0$.
\end{lem}

\begin{rem} \label{reduct1} In Theorem~\ref{caseC} we will determine
a canonical basic set for $H_n^0$. If the hypotheses of Theorem~\ref{gero}
concerning Lusztig's conjectures were known to hold in general for type
$B_n$, then Lemma~\ref{reduct0} gives a canonical basic set for $H_n$.
\end{rem}

\begin{thm} \label{caseC} Recall that $H_n^0=H_{\C}(W_n,Q_0,q_0)$ where 
$Q_0=\zeta_l^b$ and $q_0=\zeta_l^a \neq 1$ are such that $Q_0=-q_0^d$ 
for some $d \in \Z$. Let $e\geq 2$ be the multiplicative order of $q_0$
and let $p_0\in \Z$ be such that
\[ d+p_0e<\frac{b}{a}< d+(p_0+1) e.\]
(Note that the above conditions imply that $b/a \not\equiv d \bmod e$.) 
Then the set
\[ \cB^0=\Phi_{e,n}^{(d+p_0e,0)}\]
is a canonical basic set for $H_n^0$ where $\Phi_{e,n}^{(d+p_0e,0)}$ is defined
in Definition~\ref{defuglov}. The required map $\beta \colon \Lambda_n^\bu
\rightarrow \Pi_n$ such that $\cB^0=\beta(\Lambda_n^\bu)$ is given by the
map $\kappa$ in Corollary~\ref{otherui}.
\end{thm}

\begin{proof} We can identify $H_n^0$ with an Ariki--Koike algebra as in 
Remark~\ref{ident}. Hence, by  Corollary~\ref{otherui}, the decomposition 
matrix of $H_n^0$ is given by the specialisation at $q=1$ of the canonical 
basis for the highest weight module $\cM^{\bu}$ where $u_1=d+p_0 e$ and 
$u_2=0$. Now, under the isomorphism $H_n^0 \cong H_{2,n}^{\bu}$, the Specht 
module for $H_n^0$ labelled by a bipartition $\ulambda$ is isomorphic to 
the Specht module for $H_{2,n}^{\bu}$ labelled by $\ulambda$.

We will now use the same strategy as in \cite{Jac4} to prove the theorem.
We must show that for all $\umu\in \Phi_{e,n}^{(d+p_0e,0)}$:
\begin{equation*}
G(\umu,\bu)=|\umu,\bu\rangle+\sum_{\atop{\ulambda\in \Pi_n^r}{\ba_\ulambda>
\ba_\umu}} d_{\ulambda,\umu}(q)\, |\ulambda,\bu\rangle,\tag{$*$}
\end{equation*}
Again, to prove ($*$), it is sufficient to show that the matrix of the
involution on the Fock space $\cF^{\bu}$ is lower unitriangular with
respect to the $\ba$-value (see \cite[Theorem 4.6]{Jac4}). Hence, we want
to show that for all $\ulambda\in{\Pi_n^2}$, we have:
\begin{equation*}
\overline{|\ulambda,\bu\rangle}=|\ulambda,\bu\rangle+\text{sum of }
|\umu,\bu\rangle \text{ with } \ba_\ulambda<\ba_\umu.\tag{$**$}
\end{equation*}

The proof of ($**$) is rather long, but it is entirely analogous to the proof
of \cite[Theorem 4.6]{Jac4}. We only give the main arguments needed in this
proof.

First, note that the formula in Section \ref{sec:afunc} shows how we can
compute the values $\ba_\ulambda$.  Put $m^{(1)}=b/a$ and $m^{(2)}=0$. Let
$\ulambda\in{\Pi^2_n}$, $\umu\in{\Pi^2_{n+1}}$ and $\unu\in{\Pi^2_{n+1}}$ and
assume that there exists nodes $\gamma_1=(a_1,b_1,c_1)$ and
$\gamma_2=(a_2,b_2,c_2)$ such that
\[ [\umu]=[\ulambda] \cup \{\gamma_1\} \quad{\text{and}}\quad  [\unu]=[\ulambda] \cup
\{\gamma_2\}.\] Assume in addition that we have:
\[\lambda^{(c_1)}_{a_1}-a_1+m^{(c_1)}>\lambda^{(c_2)}_{a_2}-a_2+m^{(c_2)}.\]
Then it is easy to see that $\ba_\unu>\ba_\umu$ (we have a similar property
when $a$ divides $b$ in \cite[Proposition 4.3]{Jac4}).

Now,  Let $\ulambda \in{\Pi_n^2}$. Then the decomposition of
$\overline{|\ulambda,\bu\rangle}$ as a linear combination of
${|\umu,\bu\rangle}$ with $\umu\in{\Pi_n^2}$ can be obtained by using certain
rules defined by Uglov \cite[Proposition 3.16]{Ug}. These rules show that a
bipartition $|\umu,\bu\rangle$ appearing in the decomposition of
$\overline{|\ulambda,\bu\rangle}$ is obtained from another bipartition
$|\unu,\bu\rangle$, which is known by induction, by removing  a ribbon $R$ in
$\umu$ and adding a ribbon $R'$ of same size in the resulting bipartition. We
want to show that $\ba_\nu<\ba_\mu$ and the result will follow by induction.
Assume that the foot (that is the bottom-right most square) of $R$ is on the
part $\nu_{j_1}^{(i_1)}$ and that the foot of the ribbon $R'$ is on the part
$\mu_{j_2}^{(i_2)}$. Then the key property is the following:
\begin{itemize}
\item if $i_1=2$ and $i_2=1$, we have $\nu^{(i_1)}_{j_1}-j_1\geq
\nu^{(i_2)}_{j_2}-j_2+d+(p_0+1)e,$
\item if $i_1=1$ and $i_2=2$, we have $\nu^{(i_1)}_{j_1}-j_1+d+p_0 e\geq
\nu^{(i_2)}_{j_2}-j_2,$
\item if $i_1=i_2$, we have $\nu^{(i_1)}_{j_1}-j_1>\nu^{(i_2)}_{j_2}-j_2.$
\end{itemize}
The proof of this property is obtained by studying Uglov rules \cite{Ug} and is
analogous to the proof of \cite[Lemma 4.5]{Jac4}. Now, since
\[ 0<\frac{b}{a}-(d+p_0 e)<e,\]
we obtain:
\[\nu^{(i_1)}_{j_1}-j_1+m^{(i_1)}>\nu^{(i_2)}_{j_2}-j_2+m^{(i_2)}.\]
Then, we can conclude by induction exactly as in \cite[Section 4.B]{Jac4} .
\end{proof}

\begin{exmp} \label{eqparam} Assume that we have $a=1$. Then $l=e\geq 2$
must be an even number and we can take $d=b+e/2$. Then we have:
\[ b+\frac{e}{2}-e<\frac{b}{a}<b+\frac{e}{2}\]
and so $p_0=-1$. Hence,  by Theorem~\ref{caseC} and Remark~\ref{invar},
the set
\[ \cB^0=\Phi_{e,n}^{(b-e/2,0)}=\Phi_{e,n}^{(b,e/2)}\]
is a canonical basic set for $H_n^0$.

In particular, in the case $b=1$ (the ``equal parameter case''), we have
$\cB^0=\Phi_{e,n}^{(1,e/2)}$, the set of ``FLOTW bipartitions'' as in
Remark~\ref{kleshdef}(a). Thus, we recover the result shown in \cite{Jac2}.
If $b=0$, we obtain $\cB^0=\Phi_{e,n}^{(0,e/2)}$, and we recover the
result shown in \cite{Jac1}.
\end{exmp}

\begin{exmp} \label{bodda2} Assume that $a=2$ and that there exists a
nonnegative integer $r$ such that $b=2r+1$. Then $l$ must be even and
$e=l/2$ is an odd number. Then we have $d\equiv r+(1-e)/{2} \bmod e$.
Now, since
\[ r+(1-e)/2<r+1/2<r+(1+e)/2,\]
Theorem \ref{caseC} implies that $\cB^0=\Phi_{e,n}^{(r+(1-e)/2,0)}$
is a canonical basic set for $H_n^0$.
\end{exmp}

\begin{exmp} \label{bigasympt} Assume that we have $b>a(n-1)+e$.
Then we have: $d+p_0 e>n-1$. Hence  the canonical basic set
$\cB^0=\Phi_{e,n}^{(b-e/2,0)}$ for $H_n^0$ coincides with the set of Kleshchev
bipartitions by Remark~\ref{kleshdef}. Note that in this case, we have
$b>a(n-1)$, which was already dealt with in Theorem~\ref{asymp}.
\end{exmp}

\begin{rem} Following \cite{Jac4}, it would be possible to state a version 
of Theorem~\ref{caseC} which is valid for an Ariki--Koike algebra 
$H_{n,\zeta}^{\bu}$ as in Definition~\ref{arko} where $r\geq 3$. The 
$\ba$-invariants in this case are derived from the Schur elements as 
computed by Geck--Iancu--Malle \cite{GIM}. We omit further details.
\end{rem} 


\end{document}